\documentclass[letterpaper]{article}
\usepackage{amsmath}
\usepackage{amsthm}
\usepackage{amssymb}
\usepackage{graphicx}
\usepackage{tikz}

\newtheorem{corollary}{Corllary}

\theoremstyle{remark}

\title{Rigid isotopy classification of real degree~$4$ planar rational curves with only real nodes\\(An elementary approach)}
\author{Shane D'Mello}
\begin{document}
\maketitle
\begin{abstract}
  We prove an elementary method to classify, up to rigid isotopy, all nodal degree~$4$ real 
rational curves in $\mathbb{RP}^2$ that have only real double points. We
show how to associate a chord diagram to a nodal real degree~$4$ planar rational curve and prove that such curves are defined up to rigid isotopy by their chord diagrams. 
\end{abstract}

\section{Introduction}\label{s1}
All kinds of topological classification problems on real 
algebraic plane projective curves of degree 4 have been solved, see
D.A. Gudkov's survey paper ~\cite{gudkov}. However, the complexity of the 
results and just the number of inequivalent curves (117 types
of irreducible curves), is a motivation to look for parts of the
classification that can help in understanding the entire
picture. 

Algebraic curves that have only the simplest singularities, i.e. non-degenerate
double points, are called {\it nodal curves.\/} 
Rational real nodal curves occupy a very special place 
in the space of all irreducible real algebraic plane projective curves of 
degree 4: nodal curves of higher genus can be obtained from
them by small perturbations, and curves with more complicated
singularities can be obtained from nodal curves by degenerations. 
Both perturbations and degenerations for the curves under consideration 
are easy to see, although the complete list of the results is huge. 

Topological classifications of nodal rational real algebraic plane
projective curves of degree 4 were known, but the result was formulated as
just a list of pictures. 

The results of this paper provide a key to understanding the list.
We prove that the rigid isotopy type of a nodal rational 
real plane projective curve with only real nodes is determined 
by the chord diagram which describes the double points and that every 
chord diagram with at most 3 chords is realized by such curve.
This gives a new elementary proof for the rigid isotopy classification. It is straightforward to extend this classification to \textit{all} real nodal rational degree~$4$ curves, but the details of this extension will be added in a forthcoming version of this paper.

\subsection{Plane projective rational curves}\label{s1.1}
A plane projective curve $C$ is called a rational curve  
if it can be presented as the image of a regular map of the projective
line. Such a map $\theta$ is defined by formulae 
\begin{equation}\label{eq1}
[s:t]\mapsto [p_0(s,t):p_1(s,t):p_2(s,t)]
\end{equation} 
where $p_0,p_1,\mathrm{and\ }p_2$ are homogeneous polynomials of the same  
degree~$d$ with no common root. In this paper we study such curves 
of degree $d=4$ defined over the field $\mathbb R$ of real numbers.

The set of real points of $C$ is denoted by $\mathbb{R}C$ 
and the set of complex points of $C$ by $\mathbb{C}C$.
These sets are the images of maps 
$\mathbb R\theta:\mathbb{RP}^1\to\mathbb{R}C\subset\mathbb{RP}^2$ and 
$\mathbb C\theta:\mathbb{CP}^1\to\mathbb{C}C\subset\mathbb{CP}^2$ defined 
by \eqref{eq1}. 

\subsection{Real nodes}\label{s1.2}
If the preimage $\mathbb C\theta^{-1}(x)\subset\mathbb{CP}^1$ of a  
point $x\in\mathbb{C}C$ consists of two points and the
corresponding branches of $\mathbb{C}C$ are transversal, then that point is 
called a {\it node.} 

A real node $x\in\mathbb RC$ is said to be {\it non-solitary\/} if
$\mathbb C\theta^{-1}(x)\subset\mathbb CP^1$ consists of two real 
points. Then $x$ is an 
intersection point of two branches of $\mathbb RC$ transversal to 
each other in $\mathbb RP^2$.

A real node $x\in\mathbb RC$ is said to be {\it solitary\/} if 
$\mathbb C\theta^{-1}(x)\subset\mathbb CP^1$ consists of two imaginary points. Then $x$ is an intersection
point of two complex conjugate branches of $\mathbb CC$ transversal to each
other in $\mathbb{CP}^2$, but as a point of $\mathbb RC$ this point is   
isolated. 

  A generic rational complex plane projective curve of degree~$d$ has only 
nodes as singularities and they are $\frac{(d-1)(d-2)}{2}$ in number. 
Then its parametrizations $\mathbb R\theta$ and $\mathbb C\theta$ are 
generic immersions of the circle $\mathbb{RP}^1$ to $\mathbb{RP}^2$ and of
the 2-sphere $\mathbb{CP}^1$ to $\mathbb{CP}^2$, respectively. 

In particular, a generic rational complex plane projective curve of degree 
$4$ has $3$ nodes. We will consider those degree $4$ real rational plane
projective curves whose nodes are all real. Each of them may be solitary 
or non-solitary.

\subsection{Chord diagrams}\label{s1.3}
For a generic rational real plane projective curve $C$, the position of real 
non-solitary nodes on $\mathbb RC$ can be encoded into
a graph formed by a circle and a few chords. The circle is a copy of
$\mathbb{RP}^1$ parametrizing $C$ and each chord represents a non-solitary
real node of $C$ connecting the corresponding pair of points on the circle. 
The whole graph is called the {\it real chord diagram\/} of $C$.  

In the presence of imaginary nodes, one can enlarge the real chord diagram by
incorporating information about the imaginary nodes. The circle is
extended to a 2-sphere in which the original circle is the equator,
a real solitary node corresponds to a chord with end-points symmetric with
respect to the equatorial plane, an imaginary node corresponds to a chord
connecting non-symmetric points that do not belong the equator. 
However, as the curve is real, the set of
all end-points of all the chords is invariant under the symmetry in the 
equatorial plane. A chord with non-symmetric end-points can be of one of
two types: it connects either points of the same hemi-sphere or of the
different hemi-spheres. However, in this version of the paper we will restrict our attention to only those extended chord diagrams whose chords connect symmetric end points i.e. chord diagrams that correspond to curves with only real nodes. 

A {\it rigid isotopy\/} of a generic real plane projective curve is a
continuous one-parameter deformation $C_t$ with $t\in[0,1]$ consisting of
generic real plane projective curves of the same degree. Rigid isotopy is
an equivalence relation. A real chord diagram is preserved under a rigid 
isotopy.

There are two major problems related to using chord diagrams as rigid 
isotopy invariants. Firstly, it is not known which chord diagrams are 
realizable by curves of degree $d$. Secondly, it is not true that two
generic real plane projective curves that have the same chord diagram are
rigid isotopic.    

Below we solve both problems for curves of degree 4. 

\subsection{The main result}\label{s1.4} {\it Any real chord diagram 
with at most 3 chords is realizable by a nodal rational real plane 
projective curve of degree 4. 
Two nodal rational real plane projective
curves of degree 4 with only real nodes are rigidly isotopic if and only if they 
have the same real chord diagram.}
\medskip

\begin{proof}
Realize a given chord diagram with at most three chords, on the projective
line $\mathbb{RP}^1$. If the number of chords in it is less than three, then
complete it by chords connecting complex conjugate points on the sphere 
$\mathbb{CP}^1$. Let the pairs of end points of the chords be 
$[\alpha_i:\beta_i]$ and $[\gamma_i:\delta_i]$ for $i=0,1,2$.
Then the polynomials 
\begin{equation}
  \label{param}
  \begin{aligned}
p_0(s,t)=&((\beta_1s-\alpha_1t)(\delta_1s-\gamma_1t)(\beta_2s-\alpha_2t)
(\delta_2s-\gamma_2t))\\
p_1(s,t)=&((\beta_0s-\alpha_0t)(\delta_0s-\gamma_0t)(\beta_2s-\alpha_2t)
(\delta_2s-\gamma_2t))\\
p_2(s,t)=&((\beta_0s-\alpha_0t)(\delta_0s-\gamma_0t)(\beta_1s-\alpha_1t)
(\delta_1s-\gamma_1t))
\end{aligned}
\end{equation}
  define a parametrization of a generic real plane projective curve of degree
4 with the nodes at the points $z_0=[1:0:0]$, $z_1=[0:1:0]$, and $z_2=[0:0:1]$
and the prescribed chord diagram. So any chord diagram is realizable as the chord diagram corresponding to rational curve of degree~$4$. 

We now show that chord diagrams are a strong rigid isotopy invariant for curves of degree~$4$. Let $C$ be a generic real plane projective curve with three real 
nodes. The nodes cannot be co-linear, because then a line containing them 
would have an intersection number with $C$ of at least 6, which would contradict 
Bezout's theorem. Therefore, there is a projective transformation that can transform the curve so that the nodes of $C$ are the points $z_0=[1:0:0]$, $z_1=[0:1:0]$, and $z_2=[0:0:1]$.
  
  Let $[\alpha_i:\beta_i]$ and $[\gamma_i:\delta_i]$ be in the preimage of $z_i$. Let $q_i$ be the quadratic polynomials defined by 
\begin{equation}
  \begin{aligned}
  q_0(s,t)=&(\beta_0s-\alpha_0t)(\delta_0s-\gamma_0t) \\
  q_1(s,t)=&(\beta_1s-\alpha_1t)(\delta_1s-\gamma_1t) \\
  q_2(s,t)=&(\beta_2s-\alpha_2t)(\delta_2s-\gamma_2t) 
\end{aligned}
\end{equation}
If the points $[\alpha_i:\beta_i]$ and $[\gamma_i:\delta_i]$ are imaginary 
then they must occur in  conjugate pairs, so each $q_i$ is a real polynomial.

  The main observation is that when we consider $\theta$ explicitly as $[p_0:p_1:p_2]$, where $p_i$ are homogeneous polynomials of degree $4$, then the fact that $[1:0:0]$ is a node means that $q_0$ divides $p_1$ and $p_2$. Similarly, $q_1$ divides $p_0$ and $p_2$, and $q_2$ divides $p_0$ and $p_1$. 
  This proves that $p_0=c_0q_1q_2$, $p_1=c_1q_0q_2$, and $p_2=c_2q_0q_1$ where each $q_i$ is quadratic polynomials and each $c_i$ is a constant, which can all be made $1$ by a projective transformation. 
  
  Therefore, by a projective transformation, the curve has been transformed to one whose parametrization has coordinates $p_0$, $p_1$, and $p_2$ such that $p_0=q_1q_2$, $p_1=q_0q_2$, and $p_2=q_0q_1$, where each $q_i$ is determined by the preimages of the nodes as described above. This is defined explicitly by equations~\ref{param}, where the pair $\{[\alpha_i:\beta_i],[\gamma_i:\delta_i]\}$ is the preimage of each node $z_i$.
This new curve is projectively equivalent to the original curve. What we have proved is that \textit{there is a unique generic degree~$4$ rational parametrization (up to projective transformations) with the given nodes and their preimages}. Therefore, to show that two curves are isotopic, we only need to keep a track of their nodes and the pre-images of their nodes.

Consider two curves $C_1$ and $C_2$ with the same chord diagram.  They therefore have the same number of solitary nodes. Denote their explicit parametrizations by $[p_0^1:p_1^1:p_2^1]$ and $[p_0^2:p_1^2:p_2^2]$ respectively. Use a projective transformation to transform each of the curves in such a way that their nodes are sent to $[1:0:0]$, $[0:1:0]$, and $[0:0:1]$, and the solitary nodes of both curves are sent to the same points. By what we have discussed above, their resulting coordinates are of the form $[q_1^iq_2^i:q_0^iq_2^i:q_0^iq_1^i]$ for $i=1,2$, where 
\begin{equation}
  \label{components}
  \begin{aligned}
  q_0^i(s,t)=&(\beta_0^is-\alpha_0^it)(\delta_0^is-\gamma_0^it) \\
  q_1^i(s,t)=&(\beta_1^is-\alpha_1^it)(\delta_1^is-\gamma_1^it) \\
  q_2^i(s,t)=&(\beta_2^is-\alpha_2^it)(\delta_2^is-\gamma_2^it) 
\end{aligned}
\end{equation}

A path $[\alpha_i^t:\beta_i^t]$ between $[\alpha_i^1:\beta_i^1]$ and $[\alpha_i^2:\beta_i^2]$ will result in a rigid isotopy as long as the $[\alpha_i^t:\beta_i^t]$ does not coincide with any other roots of any $q_i^j$.  No two of the polynomials $q_i^j$, for each $j=1,2$, can share a root, otherwise all the coordinates will vanish on that point and the curve would not be well defined. No other singularity can be formed in the process because a generic planar rational curve of degree~$4$ can have only three nodes. 

Now we will show that the since the chord diagrams are the same, we can find a path between the roots of corresponding $q_i^j$. Suppose $[\alpha_i^1:\beta_i^1]$ is imaginary, then so will $[\alpha_i^2:\beta_i^2]$. Their conjugates would then be $[\delta_i^1:\gamma_i^1]$ and $[\delta_i^2:\gamma_i^2]$ respectively. It is always possible to find a path in $\mathbb{CP}^1\setminus \mathbb{RP}^1$ between $[\alpha_i^1:\beta_i^1]$ and $[\alpha_i^2:\beta_i^2]$ that does not contain any other root of any $q_i^j$. The conjugate of this path will be a path between their conjugates $[\delta_i^1:\gamma_i^1]$ and $[\delta_i^2:\gamma_i^2]$. 

Consider the case where $[\alpha_i^1:\beta_i^1]$ is real.  Since the curves share the same chord diagram, it is also possible to find a path $[\alpha_i^t:\beta_i^t]$ between $[\alpha_i^1:\beta_i^1]$ and $[\alpha_i^2:\beta_i^2]$. For the same reason, we can find a path connecting $[\delta_i^2:\gamma_i^2]$ and $[\delta_i^2:\gamma_i^2]$, and therefore we obtain a path between  the polynomials $q_i^1$ and $q_i^2$, for each $i$. This defines a path between $p_i^1$ and $p_i^2$, which shows that the curves $C_1$ and $C_2$ are in the same rigid isotopic class because we have been able to transform them to curves that have the same nodes and pre-images of the nodes. 
 \end{proof}

 Since there are $9$ chord diagrams with at most three chords, we have:
 \begin{corollary}
   There are $9$ rigid isotopy classes of real rational curves of degree~$4$ in $\mathbb{RP}^2$ with only real nodes as singularities. 
 \end{corollary}

 Representatives of each of these rigid isotopy classes and their respective chord diagrams are listed in figure~\ref{fig:gauss}. The dots represent the positions of the solitary nodes that which will now be deduced.

\usetikzlibrary{decorations.markings}
\begin{figure}[h]
\begin{center}
  \begin{tabular}{c|c|c|c|c}
\begin{tikzpicture}[every node/.style={draw,shape=circle,minimum size=1mm,inner sep=0pt,outer sep=0pt,fill=black},scale=.5]
  \draw (0,0) circle (1cm);
  \path (0,1) node (p0) {};
  \path (.866,.5) node (p1) {};
  \path (.866,-.5) node (p2) {};
  \path (0,-1) node (p3) {};
  \path (-.866,-.5) node (p4) {};
  \path (-.866,.5) node (p5) {};
  \draw (p5) -- (p0);
  \draw (p1) -- (p2);
  \draw (p3) -- (p4);
\end{tikzpicture}&
\begin{tikzpicture}[every node/.style={draw,shape=circle,minimum size=1mm,inner sep=0pt,outer sep=0pt,fill=black},scale=.5]
  \draw (0,0) circle (1cm);
  \path (0,1) node (p0) {};
  \path (.866,.5) node (p1) {};
  \path (.866,-.5) node (p2) {};
  \path (0,-1) node (p3) {};
  \path (-.866,-.5) node (p4) {};
  \path (-.866,.5) node (p5) {};
  \draw (p5) -- (p0);
  \draw (p1) -- (p3);
  \draw (p2) -- (p4);
\end{tikzpicture}&
\begin{tikzpicture}[every node/.style={draw,shape=circle,minimum size=1mm,inner sep=0pt,outer sep=0pt,fill=black},scale=.5]
  \draw (0,0) circle (1cm);
  \path (0,1) node (p0) {};
  \path (.866,.5) node (p1) {};
  \path (.866,-.5) node (p2) {};
  \path (0,-1) node (p3) {};
  \path (-.866,-.5) node (p4) {};
  \path (-.866,.5) node (p5) {};
  \draw (p5) -- (p0);
  \draw (p4) -- (p1);
  \draw (p3) -- (p2);
\end{tikzpicture}&
\begin{tikzpicture}[every node/.style={draw,shape=circle,minimum size=1mm,inner sep=0pt,outer sep=0pt,fill=black},scale=.5]
  \draw (0,0) circle (1cm);
  \path (0,1) node (p0) {};
  \path (.866,.5) node (p1) {};
  \path (.866,-.5) node (p2) {};
  \path (0,-1) node (p3) {};
  \path (-.866,-.5) node (p4) {};
  \path (-.866,.5) node (p5) {};
  \draw (p0) -- (p3);
  \draw (p1) -- (p5);
  \draw (p4) -- (p2);
\end{tikzpicture}&
\begin{tikzpicture}[every node/.style={draw,shape=circle,minimum size=1mm,inner sep=0pt,outer sep=0pt,fill=black},scale=.5]
  \draw (0,0) circle (1cm);
  \path (0,1) node (p0) {};
  \path (.866,.5) node (p1) {};
  \path (.866,-.5) node (p2) {};
  \path (0,-1) node (p3) {};
  \path (-.866,-.5) node (p4) {};
  \path (-.866,.5) node (p5) {};
  \draw (p3) -- (p0);
  \draw (p5) -- (p2);
  \draw (p1) -- (p4);
\end{tikzpicture}\\
\hline
\includegraphics[height=1cm]{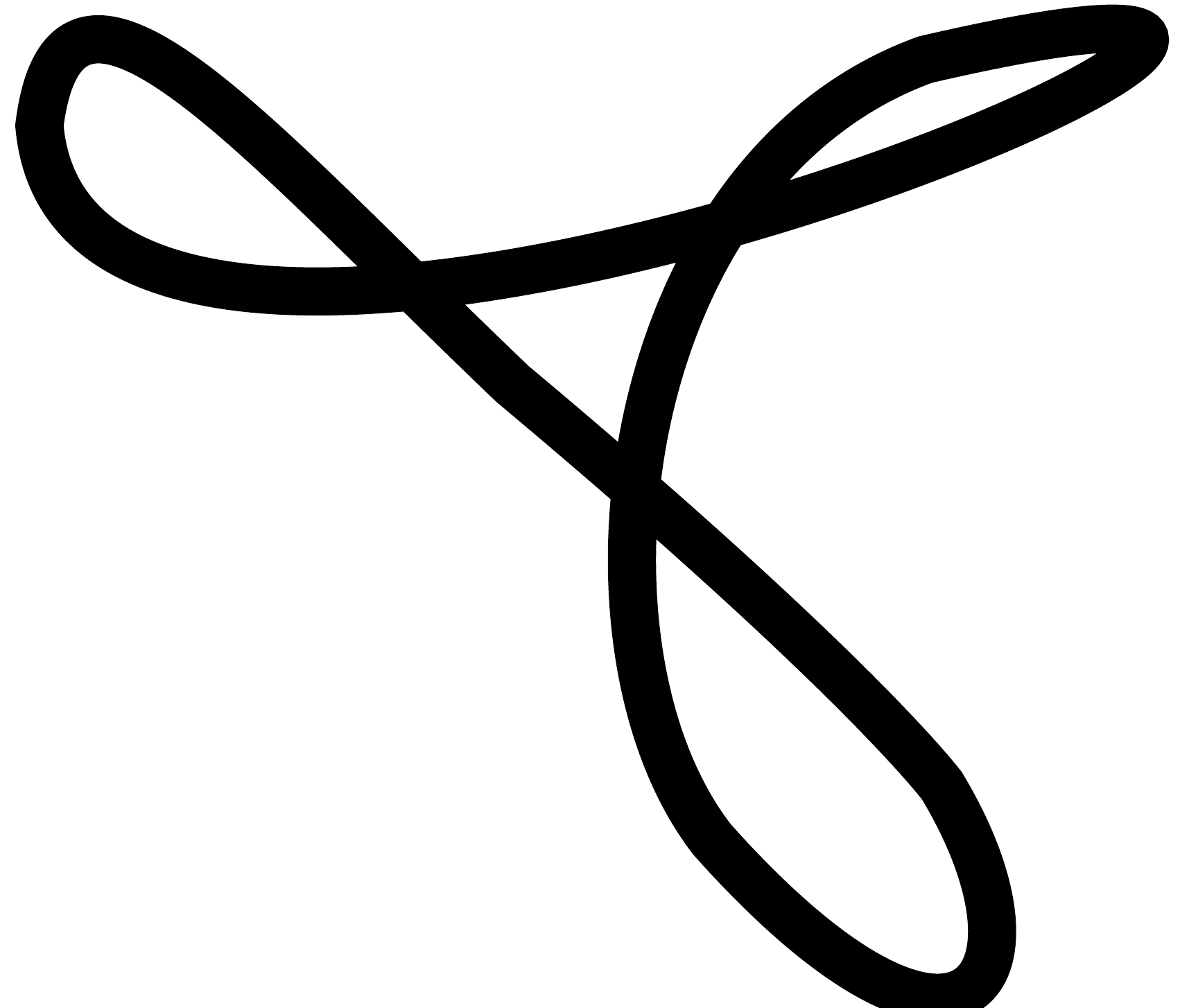}&
\includegraphics[height=1cm]{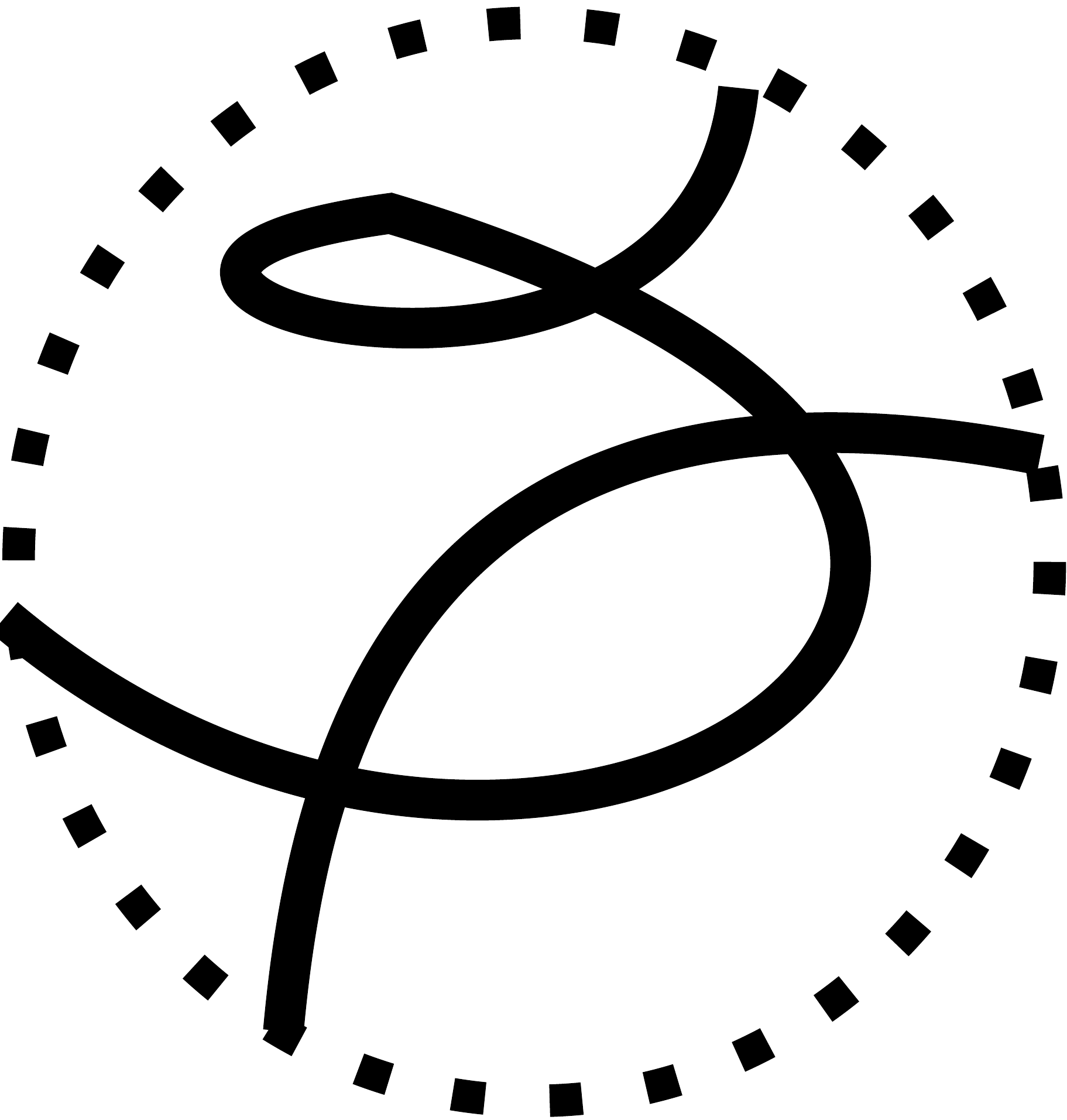}&
\includegraphics[height=1cm]{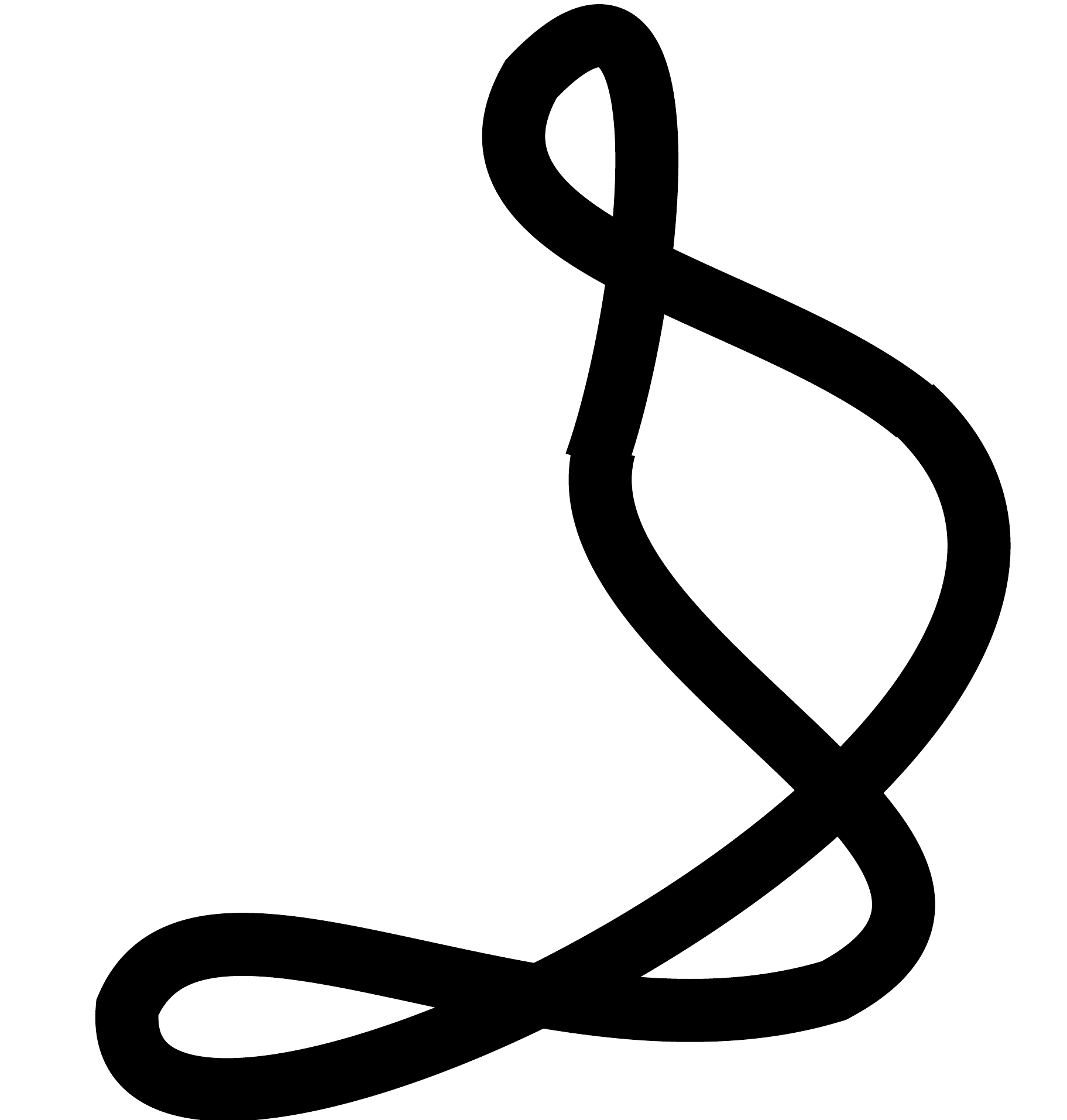}&
\includegraphics[height=1cm]{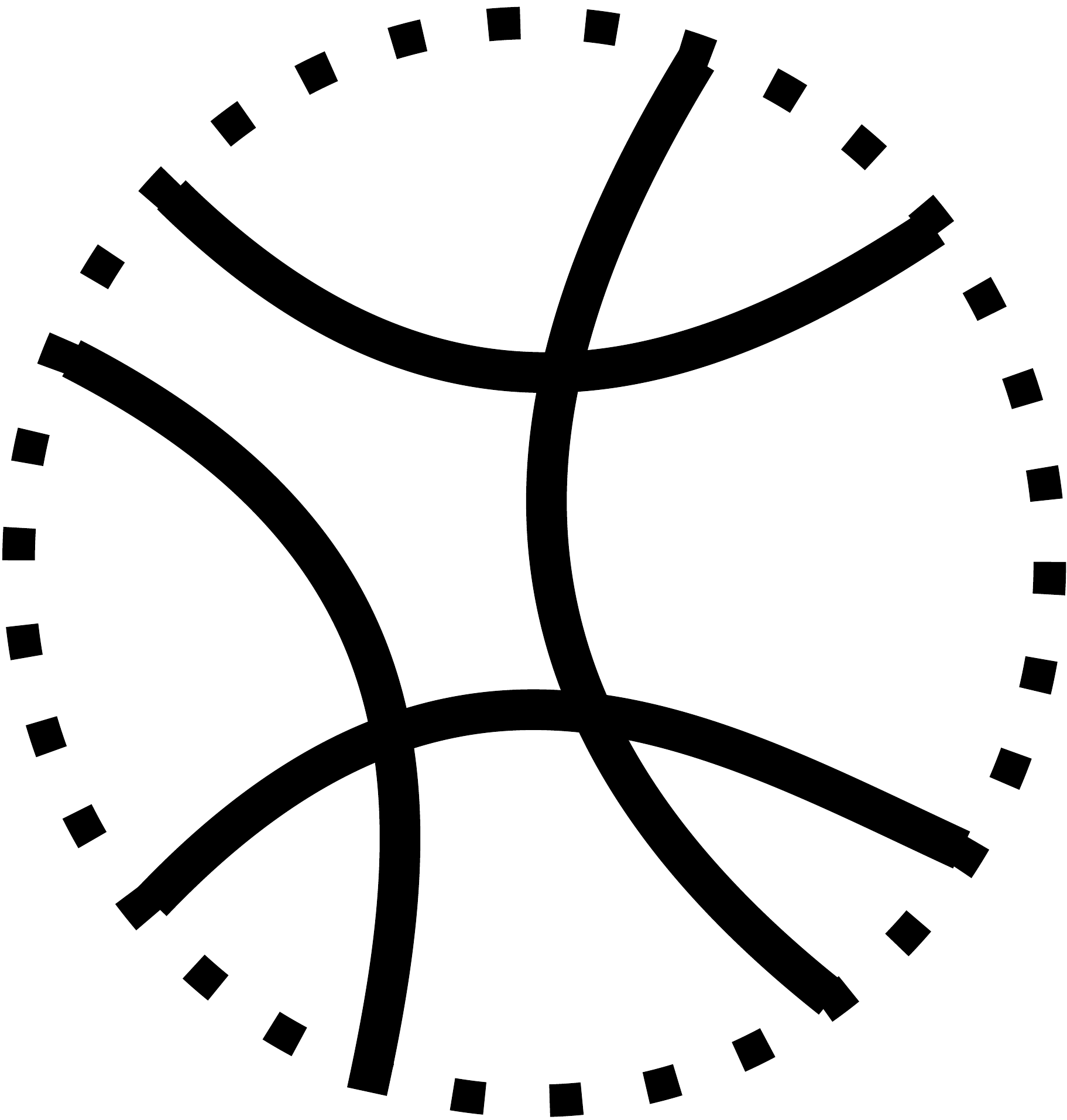}&
\includegraphics[height=1cm]{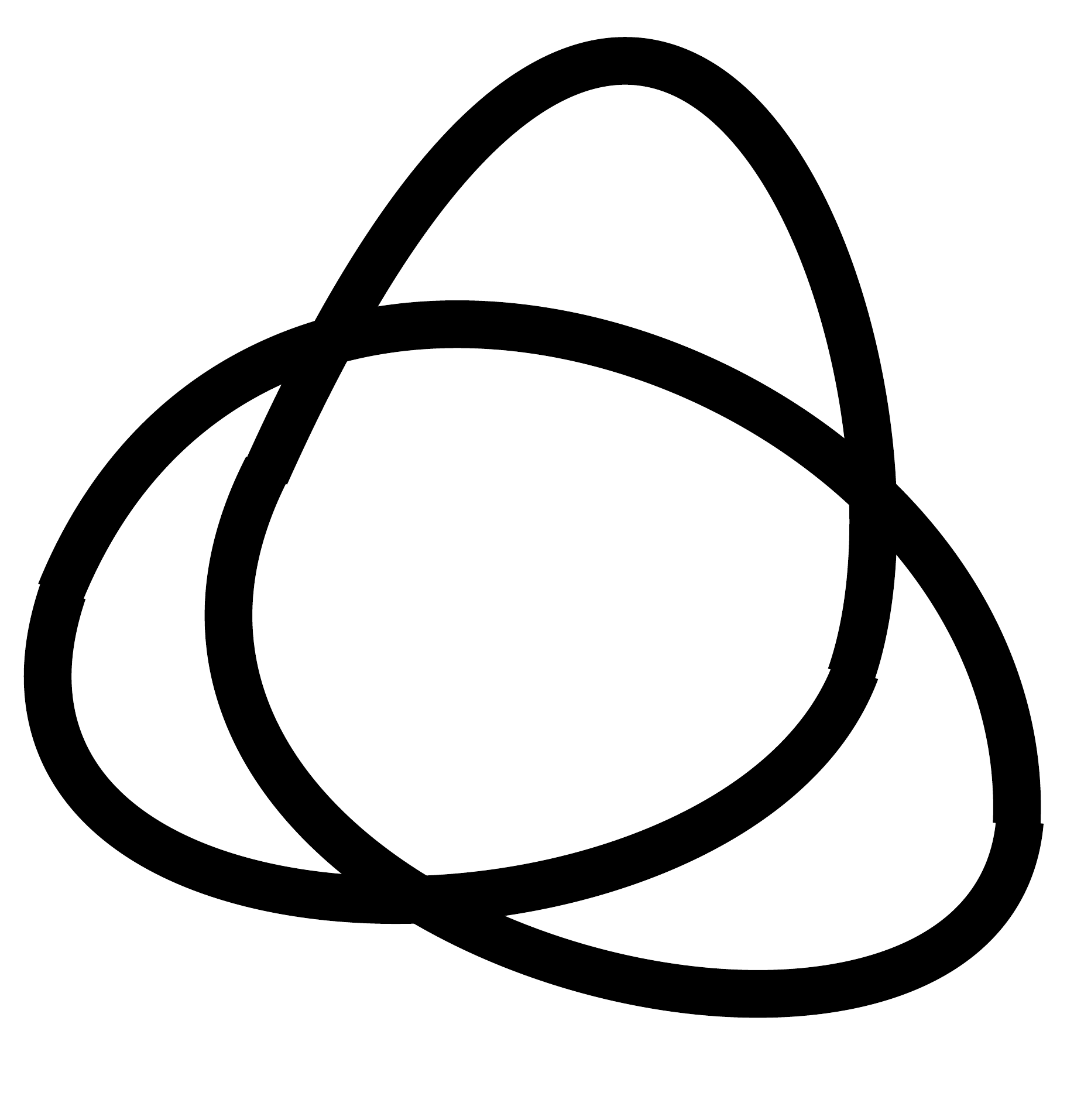}
  \end{tabular}\\
\end{center}
\begin{center}
  \begin{tabular}{c|c|c|c}
\begin{tikzpicture}[every node/.style={draw,shape=circle,minimum size=1mm,inner sep=0pt,outer sep=0pt,fill=black},scale=.5]
  \draw (0,0) circle (1cm);
  \path (0,1) node (p0) {};
  \path (1,0) node (p1) {};
  \path (0,-1) node (p2) {};
  \path (-1,0) node (p3) {};
  \draw (p0) -- (p2);
  \draw (p1) -- (p3);
\end{tikzpicture}&
\begin{tikzpicture}[every node/.style={draw,shape=circle,minimum size=1mm,inner sep=0pt,outer sep=0pt,fill=black},scale=.5]
  \draw (0,0) circle (1cm);
  \path (0,1) node (p0) {};
  \path (1,0) node (p1) {};
  \path (0,-1) node (p2) {};
  \path (-1,0) node (p3) {};
  \draw (p0) -- (p1);
  \draw (p2) -- (p3);
\end{tikzpicture}&
\begin{tikzpicture}[every node/.style={draw,shape=circle,minimum size=1mm,inner sep=0pt,outer sep=0pt,fill=black},scale=.5]
  \draw (0,0) circle (1cm);
  \path (0,1) node (p0) {};
  \path (0,-1) node (p2) {};
  \draw (p0) -- (p2);
\end{tikzpicture}&
\begin{tikzpicture}[every node/.style={draw,shape=circle,minimum size=1mm,inner sep=0pt,outer sep=0pt,fill=black},scale=.5]
  \draw (0,0) circle (1cm);
\end{tikzpicture}\\
\hline
\includegraphics[height=1cm]{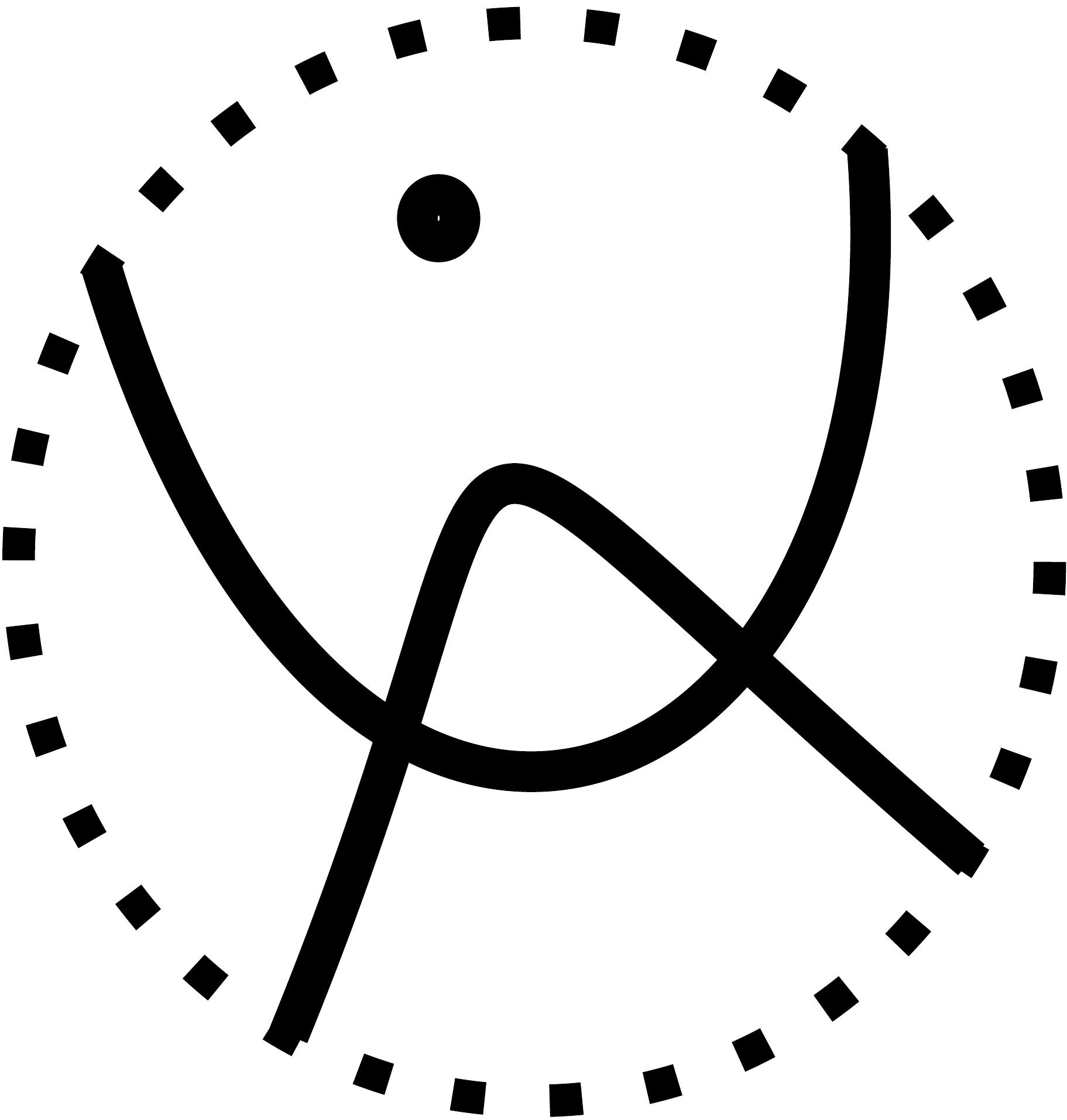}&
\includegraphics[height=1cm]{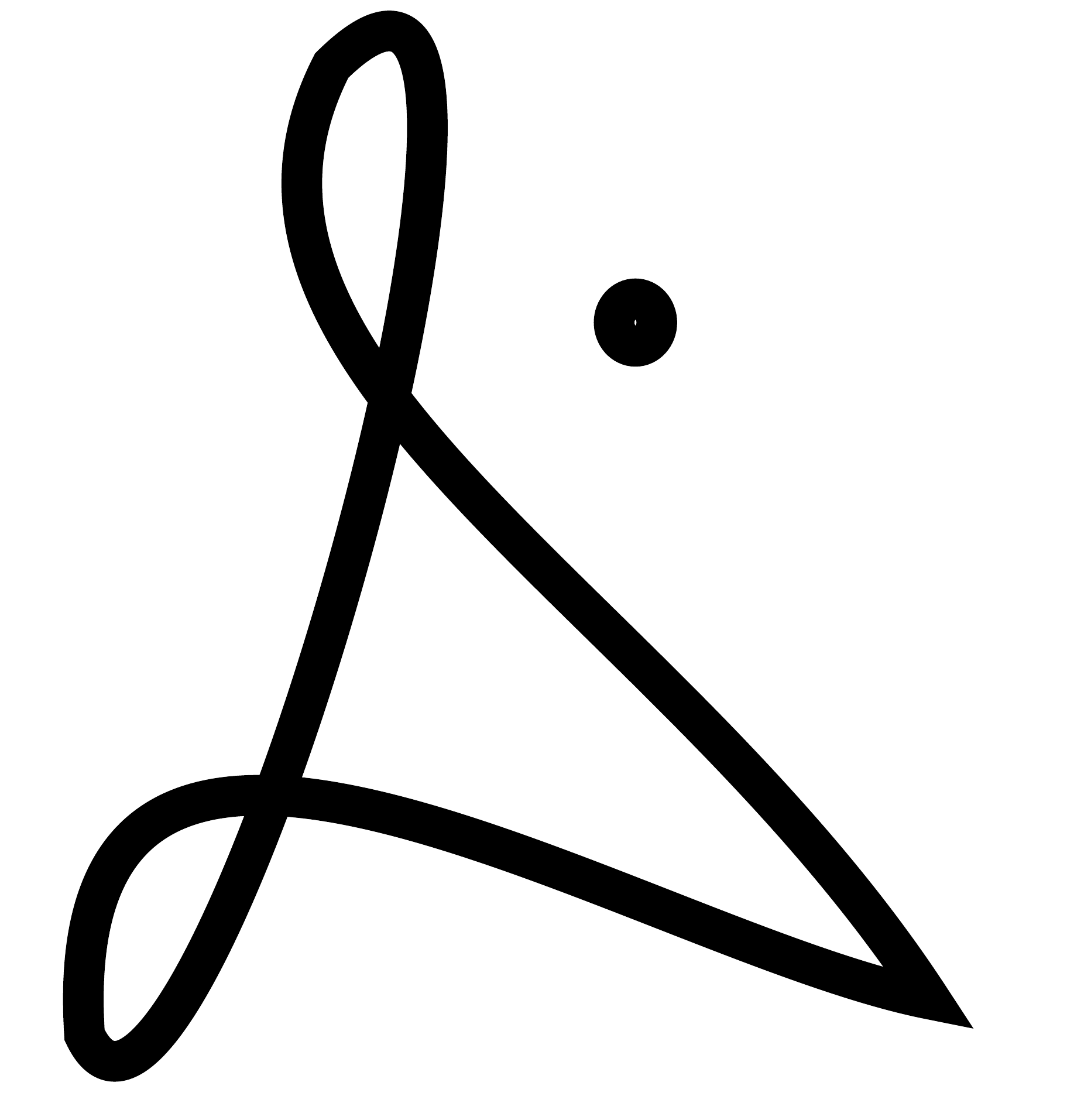}&
\includegraphics[height=1cm]{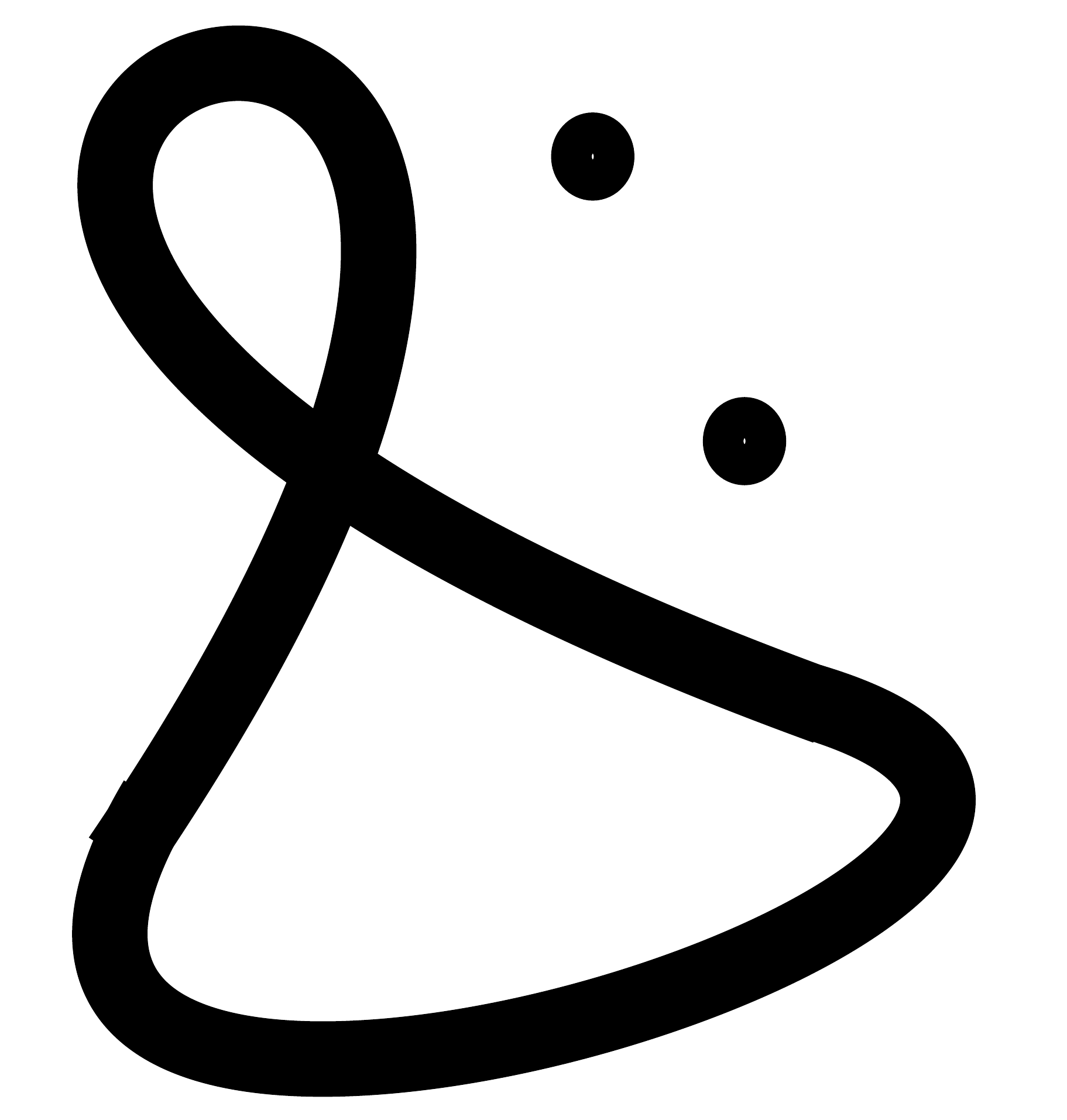}&
\includegraphics[height=1cm]{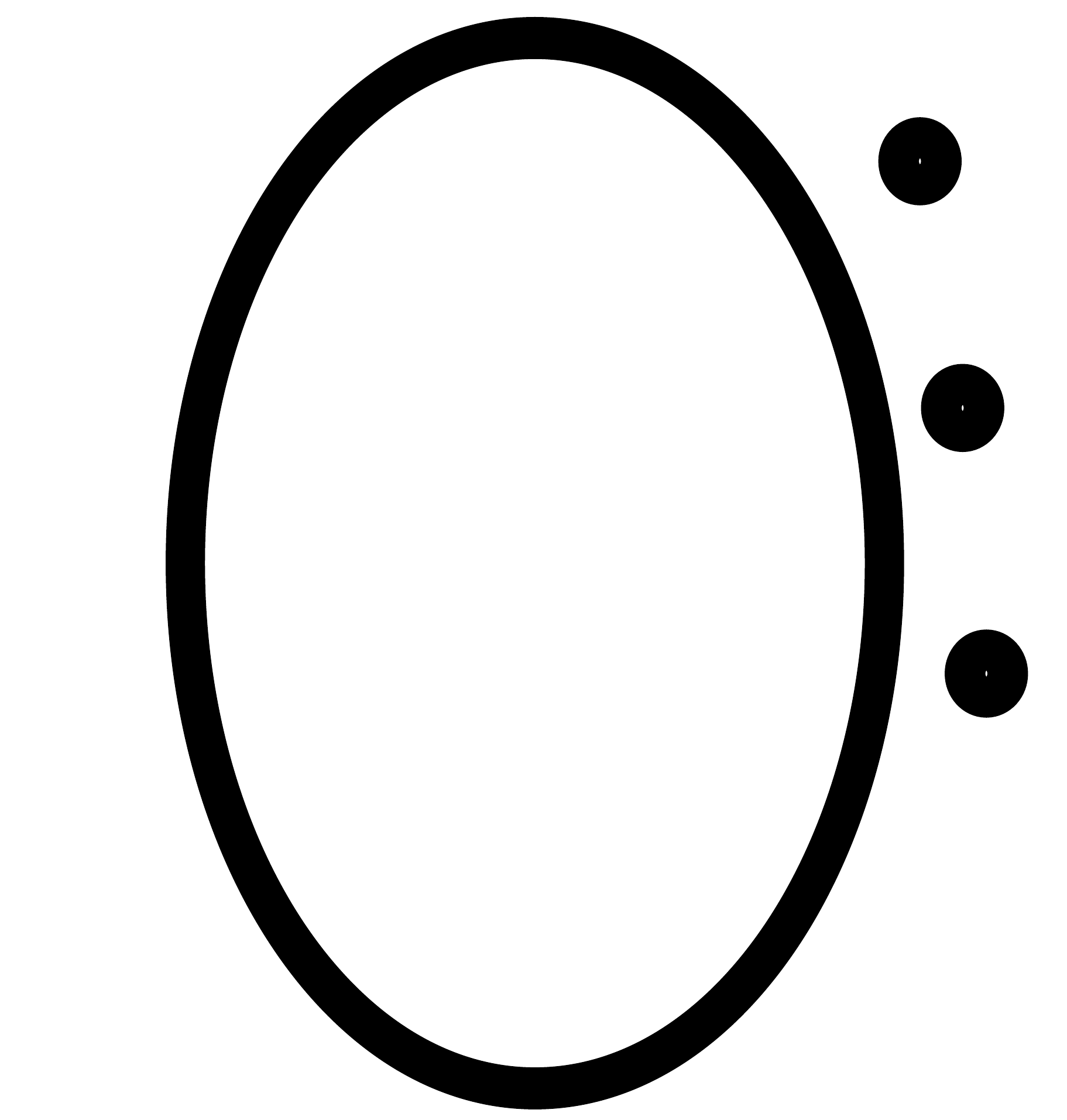}
  \end{tabular}
\end{center}
	\caption{Chord diagrams of planar curves of degree $4$ with real nodes. The lower row consists of curves with solitary nodes, whose positions are denoted by a dot.}
	\label{fig:gauss}
\end{figure}
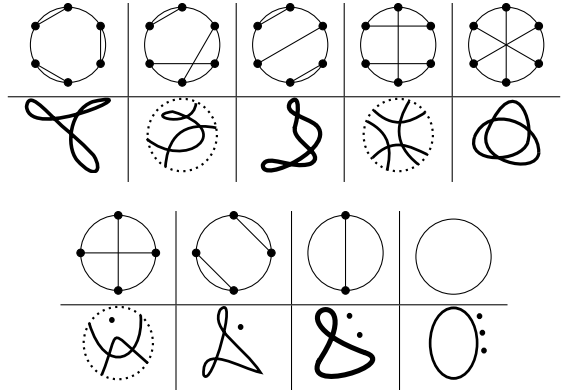

 \section{Positions of the solitary nodes}
 As we have seen, the positions of the nodes and pre-images uniquely determine the rational curve. Therefore, for a curve $C$ with a parametrization $\theta$, we should be able to uniquely determine which components of $\mathbb{RP}^3\setminus\theta(\mathbb{RP}^1)$ each of the solitary nodes must lie in. We will now determine the component and show that if the curve has more than one solitary node, then all the solitary nodes lie in the same component of $\mathbb{RP}^3\setminus\theta(\mathbb{RP}^1)$

 Consider the curves in figure~\ref{fig:gauss} that have solitary nodes. They are the last four curves of figure~\ref{fig:gauss}. Call them $C_1, C_2, C_3$, and $C_4$ and denote their respective parametrizations by $\theta_1$, $\theta_2$, $\theta_3$, and $\theta_4$. Then, $\mathbb{RP}^2\setminus \theta_2(\mathbb{RP}^1), \mathbb{RP}^2\setminus \theta_3(\mathbb{RP}^1)$, and $\mathbb{RP}^2\setminus \theta_4(\mathbb{RP}^1)$ each have only one component that is not a disk.
 
\begin{figure}[h]
  \begin{center}
    \includegraphics[height=2cm]{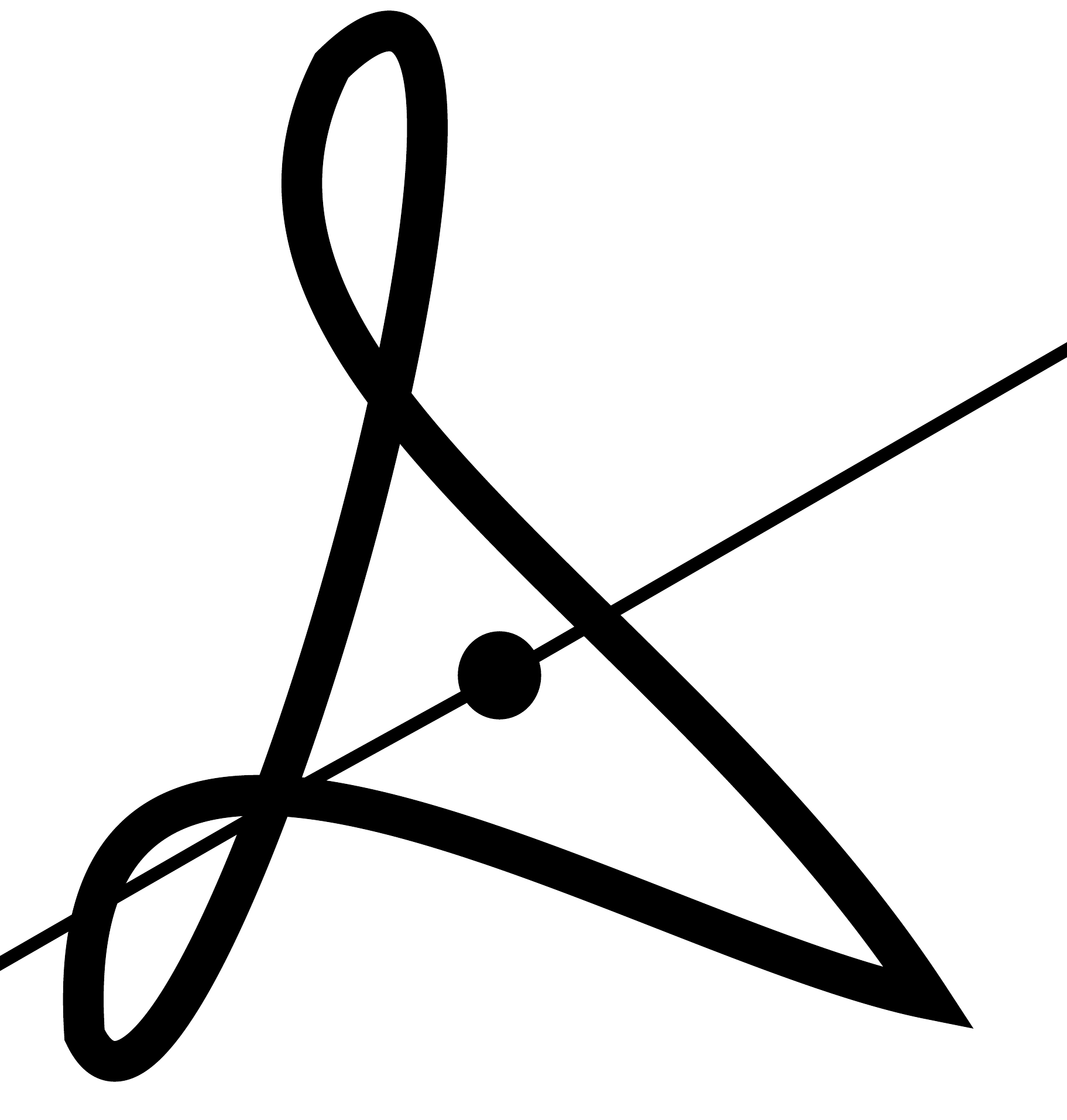}
  \end{center}
  \caption{An example of an impossible position for a solitary node because the line would intersect the curve in more than $5$ points counting multiplicity.}
  \label{fig:line}
\end{figure}

 If a solitary node lay in the interior of a disk, any line passing through it would have to intersect the boundary of the disk twice. So a line joining this solitary node and another node would intersect the curve in at least $5$ points counting multiplicity: four from the double intersections with the two nodes defining the line, and at least one more from the intersection with the boundary of the disk at a point that is not a node. So none of the solitary nodes can lie inside the disk components. Figure~\ref{fig:line} shows one example of this applied to $C_2$.

 Since the complements of $\theta_2(\mathbb{RP}^1), \theta_3(\mathbb{RP}^1)$, and $\theta_4(\mathbb{RP}^1)$, each have only one non-disk component, their solitary nodes have to lie in it as is shown for the last three curves in figure~\ref{fig:gauss}. 

We are only left with the curve $C_1$. $\theta_1(\mathbb{RP}^1)$ divides $\mathbb{RP}^2$ into three components out of which only one is a disk. The same reasoning as above will show that its node cannot lie in that disk. To deduce which of the other components it lies on, we will need to use Rokhlin's Complex Orientation Formula for non-singular curves of Type I which is briefly reviewed in the next section. For more details on Rokhlin's Complex Orientation Formula see~\cite{viro:TRAV}.

 \subsection{Review of Rokhlin's Complex Orientation Formula}
 A real algebraic curve $A$ is said to be of \textit{Type I} if $[\mathbb{R}A]=0\in\mathrm{H}_1(\mathbb{C}A)$, otherwise it is said to be of \textit{Type II}. If $A$ is of Type I, then $\mathbb{R}A$ is the boundary of each of the two halves of $\mathbb{C}C\setminus \mathbb{R}A$ which are conjugate to each other. If the halves intersect, they do so precisely at the solitary nodes.  
 
 Each half induces a natural pair of opposite orientations on their common boundary $\mathbb{R}A$. This pair of orientations is called the \textit{complex orientation} of the curve. A rational curves is an example of a Type I curve and its complex orientations coincide with the natural orientations induced on it by its parametrization. 
 
We now consider non-singular real algebraic curves defined by the zero set of a real polynomial. The real zero set $\mathbb{R}A$ of a non-singular real algebraic curve $A$ is a compact one-dimensional manifold and therefore each component is topologically a circle.  

If the complement of a component of $\mathbb{R}A$ divides $\mathbb{RP}^2$ it is called an \textit{oval}. Each oval is the boundary of a disk in $\mathbb{RP}^2$; the disk is called the interior of the oval. 
 
A pair of ovals is called \textit{injective} if one lies in the interior of the other. A pair of injective ovals bounds an annulus. If the orientation of the ovals is compatible with the orientation induced on the boundary of the annulus, it is called a \textit{positive pair}, otherwise it is called a \textit{negative pair}. The number of positive pairs is denoted by $\Pi^+$ and the number of negative pairs is denoted by $\Pi^-$. Rokhlin proved that non-singular real planar curves of Type I of degree~$d$ and with $l$ components satisfy the following formula known as \textit{Rokhlin's Compex Orientation Formula}~\cite{rokhlin}

 \begin{equation}
   2(\Pi^+-\Pi^-)=l-\frac{d^2}{4}
   \label{rokhlin}
 \end{equation}

\subsection{The curve $C_1$}
\begin{figure}[h]
  \begin{center}
    \includegraphics[height=2cm]{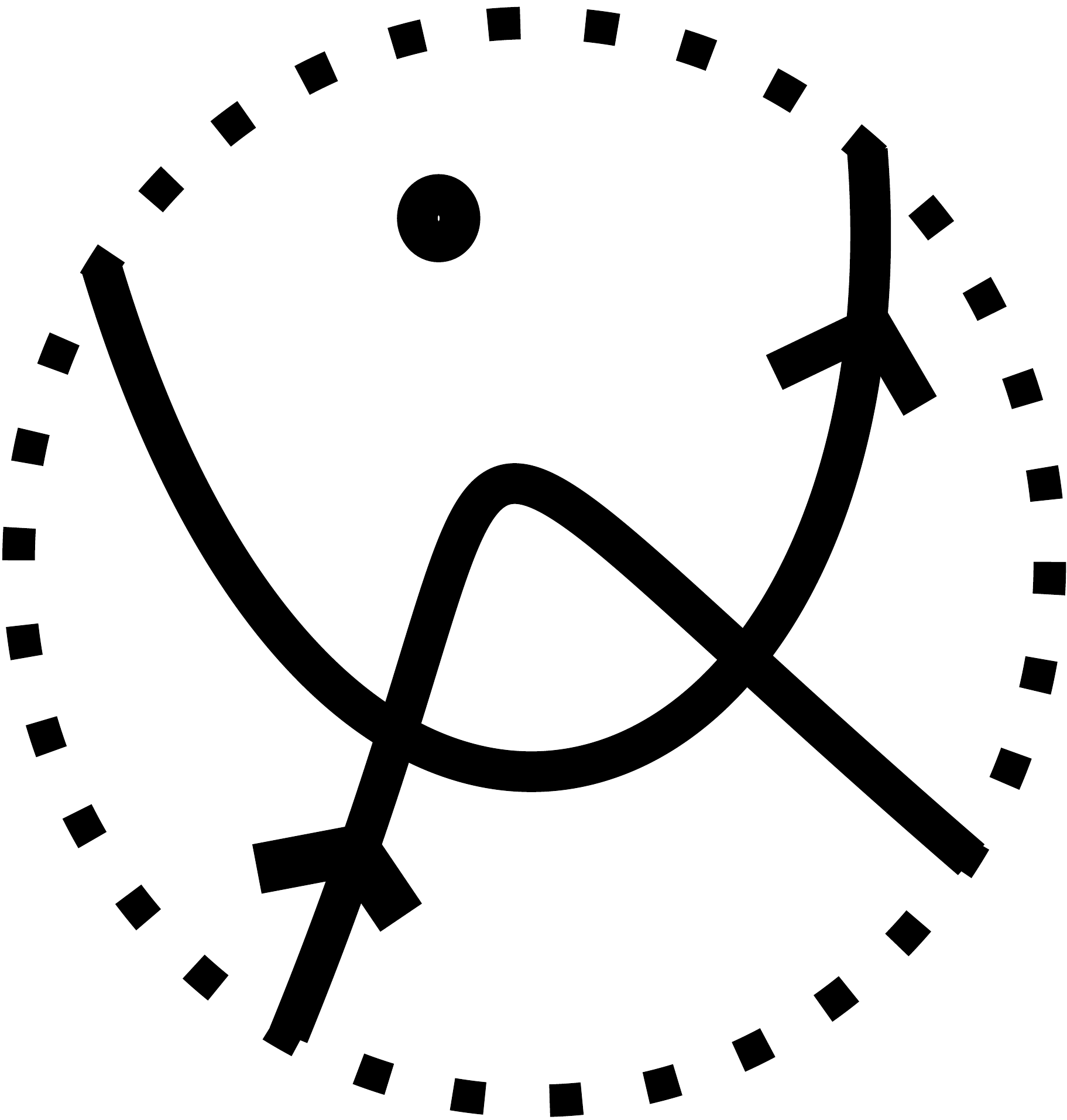}
    \includegraphics[height=2cm]{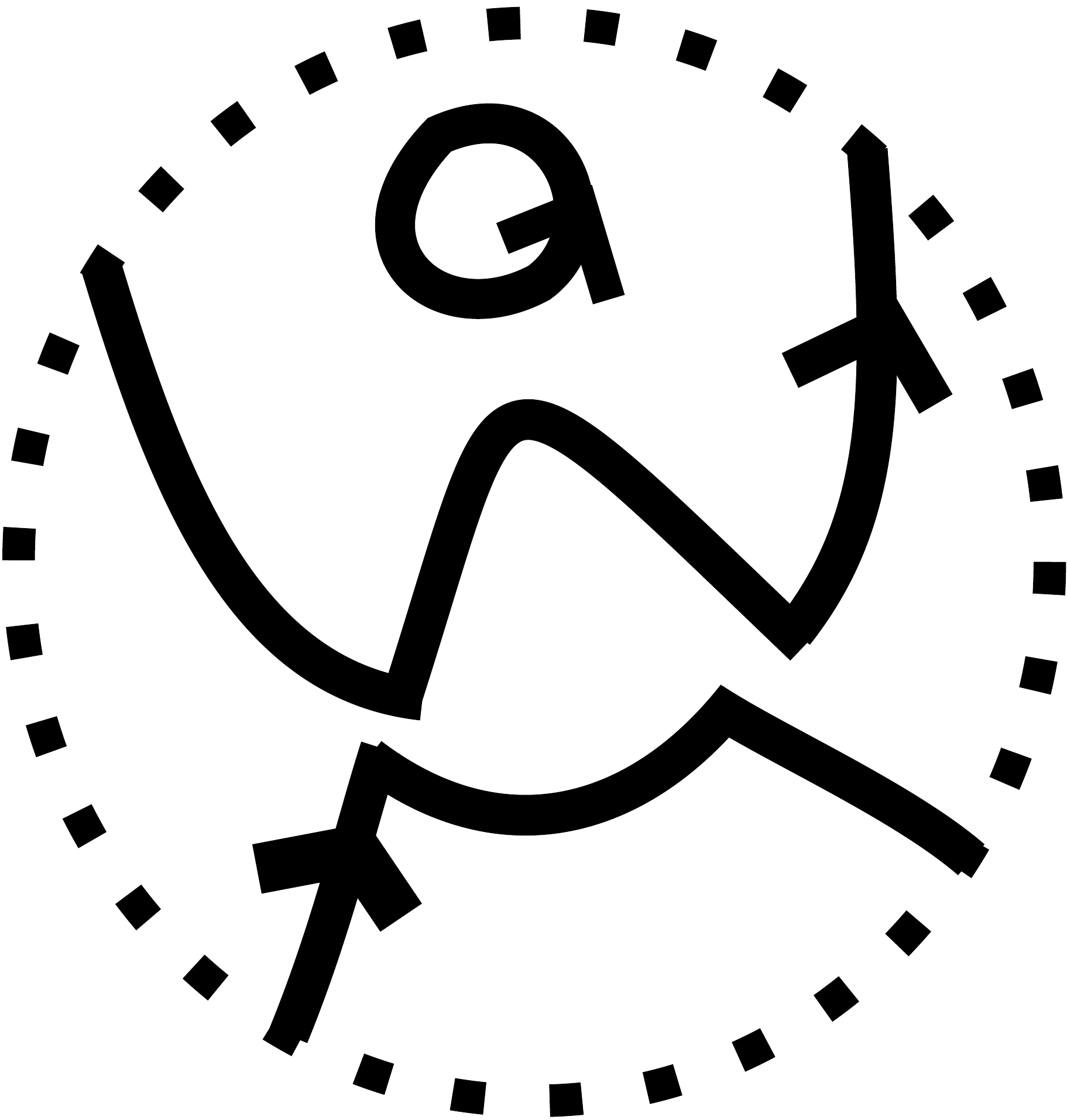}
  \end{center}
  \caption{$C_1$ and its perturbation with a complex orientation}
  \label{fig:c1}
\end{figure}

Orient $\mathbb{R}C_1$. This orientation corresponds to one of its complex orientations, which is equivalent to choosing one half of the complexification. It is possible to smoothen the real part of the curve in such a way that the double point changes to an oval and the curve continues to remain of Type I (see section 2.2 of \cite{viro:generic}).  There is only one way to resolve the crossings so that they conform to the complex orientation that will be induced. This, as shown in figure~\ref{fig:c1},  will give rise to a non-singular degree~$4$ planar curve with two ovals. One of the ovals is the result of the perturbation of the solitary node. 

If the pair of ovals is not injective, then the number of positive pairs of ovals $\Pi^+$ and the number of negative pairs of ovals $\Pi^-$ will be $0$. But the right hand side of formula~\ref{rokhlin} should be $-2$, because $d=4$ and $l=2$. There is only one possible component of  $\mathbb{RP}^2\setminus \theta_1(\mathbb{RP}^1)$ in which the solitary node can lie in, so that on perturbing the resulting ovals will be nested. The only possible position for the solitary node is shown in figures~\ref{fig:c1} and \ref{fig:gauss}. 

This completes the proof of the rigid isotopy classification of degree~$4$ projective planar rational curves with real nodes, including the positions of the isolated nodes. In a future version of this paper, this classification will be extended to \textit{all} generic rational planar curves of degree~$4$, i.e. even those that have imaginary nodes.
 
\section{Acknowledgements}
I am very grateful to Prof. Oleg Viro for providing very valuable ideas and suggestions. 
\bibliography{references}{}
\bibliographystyle{acm}
\end{document}